\documentclass[11pt]{article}
\usepackage{amsfonts,latexsym}

\renewcommand{\baselinestretch}{1.34}
\makeatletter
\def\operatorname#1{\mathop{\operator@font #1}\nolimits}%
\makeatother
\renewcommand{\a}{\mathfrak{a}}
\newcommand{\C}{\mathbb{C}}
\newcommand{\esl}{\mathfrak{sl}}
\newcommand{\g}{\mathfrak{g}}
\newcommand{\gl}{\mathfrak{gl}}
\renewcommand{\k}{\mathfrak{k}}
\renewcommand{\P}{\mathbb{P}}
\newcommand{\p}{\mathfrak{p}}
\newcommand{\R}{\mathbb{R}}
\renewcommand{\r}{\mathfrak{r}}

\newcommand{\so}{\mathfrak{so}}
\newcommand{\su}{\mathfrak{su}}
\renewcommand{\u}{\mathfrak{u}}
\newcommand{\Ad}{\operatorname{Ad}}
\newcommand{\Mat}{\operatorname{Mat}}
\newcommand{\Aut}{\operatorname{Aut}}
\newcommand{\End}{\operatorname{End}}
\newcommand{\Exp}{\operatorname{Exp}}
\renewcommand{\exp}{\operatorname{exp}}
\newcommand{\Id}{\operatorname{Id}}
\renewcommand{\Im}{\operatorname{Im}}
\renewcommand{\Re}{\operatorname{Re}}
\newcommand{\Image}{\operatorname{Image}}
\newcommand{\Ker}{\operatorname{Ker}}
\newcommand{\rank}{\operatorname{Rank}}
\newcommand{\Tr}{\operatorname{Tr}}
\newcommand{\llangle}{\langle\langle}
\newcommand{\rrangle}{\rangle\rangle}
\newcommand{\suchthat}{\mathop{\,\vert\,}}
\let\ul\underline
\let\ol\overline

\newtheorem{thm}{Theorem}

\newtheorem{cor}[thm]{Corollary}
\newtheorem{lemma}[thm]{Lemma}
\newtheorem{definition}[thm]{Definition}
\makeatletter
\newenvironment{proof}[1][]%
 {\def\proof@temp{#1}\par\noindent
  \textsc{Proof}\ifx\proof@temp\@empty\else\ (#1)\fi\hspace{1em}}
 {~~\hfill{$\Box$}\par\vspace{.2\baselineskip}}
\makeatother

\newcounter{Item}
\def\Item{\stepcounter{Item}\vspace{0.5\baselineskip}%
\noindent\makebox[2em][l]{(\theItem)}~}

\newcommand{\email}[1]{{\small\tt #1}}

\def\ftnote#1{\def\footnotemark{}\footnote{#1}\setcounter{footnote}{0}}

\title{Symmetric symplectic spaces\\ with Ricci-type 
curvature}

\author{M. Cahen\footnotemark\addtocounter{footnote}{-1}
\ftnote{For the proceedings of the 
Conference Mosh\'e Flato, Dijon, 1999.}\\
\email{mcahen\char64ulb.ac.be}\\[5pt]
\small Universit\'e Libre de Bruxelles\\
\small Campus Plaine, CP 218\\
\small bd du Triomphe\\
\small 1050 Brussels, Belgium\\[15pt]
S. Gutt\thanks{Research supported by the Communaut\'e fran\c caise de
Belgique, through an Action de Recherche Concert\'ee de la Direction de
la Recherche Scientifique.}\\
\email{sgutt\char64ulb.ac.be}\\[-8pt]
\small \hbox{\parbox[t]{1.9in}{\begin{center}%
Universit\'e Libre de Bruxelles\\[3pt]
Campus Plaine, CP 218\\[3pt]
bd du Triomphe\\[3pt]
1050 Brussels, Belgium\end{center}}}
\hbox{\parbox[t]{.6in}{\begin{center}\rm and\end{center}}}
\hbox{\parbox[t]{1.9in}{\begin{center}%
Universit\'e de Metz\\[3pt]
Ile du Saulcy\\[3pt]
57045 Metz Cedex 01, France\end{center}}}
\\[15pt]
J. Rawnsley\\
\email{j.rawnsley\char64warwick.ac.uk}\\[5pt]
\small Mathematics Institute\\
\small University of Warwick\\
\small Coventry CV4 7AL, UK}

\date{December 1999}

\begin{document}
{\renewcommand{\baselinestretch}{1} 
\maketitle\thispagestyle{empty}

\begin{abstract}
We determine the isomorphism classes of symmetric symplectic manifolds
of dimension at least 4 which are connected, simply-connected
and have a curvature tensor which has only one non-vanishing irreducible
component -- the Ricci tensor.
\end{abstract}



\vglue20pt\leftskip2pc

\parbox{30pc}{\noindent \it Mosh\'e Flato has been a close and wonderful friend and an
inspiration for us for more than twenty years. This contribution is
dedicated to him, always present in our hearts.}

}

\newpage

\section{Introduction}

On any symplectic manifold $(M,\omega)$ the space of symplectic
connections (linear connections $\nabla$ with vanishing torsion and such
that $\nabla\omega=0$) is infinite dimensional. In order to select a
smaller family of symplectic connections, a variational principle was
introduced in \cite{Bourg}. This principle has Euler-Lagrange equations
\begin{equation}\label{Intro:field}
(\nabla_Xr)(Y,Z)+(\nabla_Yr)(Z,X)+(\nabla_Zr)(X,Y)=0
\end{equation}
for all vector fields $X,Y,Z$; $r$ denotes the Ricci tensor of $\nabla$
\[
r(X,Y) = \Tr (Z\mapsto R(X,Z)Y).
\]
In \cite{Bourg} the case where $\dim M=2$ was examined in complete detail
so we shall assume throughout that the dimension is at least 4.

It was observed in \cite{CGR} that the field equations (\ref{Intro:field})
are identically satisfied if one assumes that the irreducible component
of the curvature, denoted there by $W$ (see also \cite{Vais}), vanishes
\begin{equation}\label{Intro:Weyl}
W=0.
\end{equation}
The tensor $W$ is the symplectic analogue of the Weyl or conformal
curvature of a Riemannian connection. The vanishing of $W$ (equation
(\ref{Intro:Weyl})) is equivalent to the requirement that the curvature
tensor $R$ of $\nabla$ is expressed in terms of its Ricci tensor by
\begin{eqnarray}
R(X,Y)Z &=& \frac{1}{2(n+1)}\biggl[ 2\omega(X,Y)AZ + \omega(X,Z)AY
-\omega(Y,AZ)X\nonumber\\\label{Intro:curv}
&&\qquad - \omega(Y,Z)AX + \omega(X,AZ)Y\biggr]
\end{eqnarray}
where $\dim M = 2n$, $n\ge2$, where $X,Y,Z$ are vector fields and 
where $A$ is the
Ricci tensor viewed as an endomorphism of the tangent bundle using
$\omega$:
\begin{equation}\label{Intro:endo}
r(X,Y) = \omega(X,AY).
\end{equation}
The Ricci tensor is symmetric so $A$ is an infinitesimal symplectic
endomorphism of each tangent space.
 
Equations (\ref{Intro:Weyl}) (or (\ref{Intro:curv})) imply the existence
of a 1-form $u$ on $(M,\omega)$ such that
\begin{equation}\label{Intro:formu}
(\nabla_Xr)(Y,Z) = \omega(X,Y)u(Z) + \omega(X,Z)u(Y).
\end{equation}
If $u=0$, then $\nabla r=0$ and since $R$ is expressed in terms of $r$ 
(\ref{Intro:curv}), $\nabla$ is locally symmetric.

The condition $W=0$ also appears as the integrability condition
for the almost complex structure naturally defined from a symplectic
connection on $(M,\omega)$ on the manifold ${\mathcal{J}}(M)$ of almost
complex stuctures on $M$ which are compatible with $\omega$.

In this note we prove, amongst other things, the following two results.

\begin{thm} \label{Intro:thm1}
Let $(M,\omega) = (M_1,\omega_1) \times (M_2,\omega_2)$ be symplectic
manifolds of dimension greater than zero and $\nabla = \nabla_1 +
\nabla_2$ be a symplectic connection. If $W^\nabla = 0$ then $\nabla,
\nabla_1 , \nabla_2$ are flat.
\end{thm}

\begin{thm} \label{Intro:thm2}
Let $(M,\omega,s)$ be a connected, simply-connected, symmetric
symplectic space of dimension $2n (\ge4)$; let $\nabla$ be its canonical
invariant symplectic connection and let $r$ be its Ricci tensor; let $A$
be the corresponding endomorphism
\[
\omega(X,AY) = r(X,Y).
\]
Assume $W^\nabla = 0$. Then 
\[
A^2 = \lambda \Id
\]
for some real number $\lambda$.

If $\lambda\ne0$, the transvection group $G$ of $(M,\omega,s)$ is
semisimple and, up to coverings, $M=G/K$ with either $G=SL(n+1,\R)$ and
$K=GL(n,\R)$ or $G=SU(p+1,q)$ and $K=U(p,q)$ where $\dim M=2n$, $p+q=n$.

If $\lambda=0$ and $\rank(A)>1$, the transvection group $G$ of
$(M,\omega,s)$ is neither solvable nor semisimple. The radical of $G$ is
2-step unipotent if $\rank(A)<n$ and abelian in $\rank(A)=n$. If
$\lambda=0$ and $\rank(A)=1$, the transvection group $G$ of
$(M,\omega,s)$ is solvable.
\end{thm}

\def\temprefa{\ref{Intro:thm1}}
\section{Proof of Theorem \temprefa}

Let $(M,\omega) = (M_1,\omega_1) \times (M_2,\omega_2)$ be symplectic
manifolds and $\nabla = \nabla_1 + \nabla_2$ be a symplectic connection.
Then $R(X,Y)Z = R_1(X_1,Y_1)Z_1 + R_2(X_2,Y_2)Z_2$ where $X=X_1+X_2$,
$Y=Y_1+Y_2$, $Z=Z_1+Z_2$ and suffices indicate components tangent to
$M_1$ and $M_2$, respectively. Then also $r(X,Y) = r_1(X_1,Y_1) +
r_2(X_2,Y_2)$. On the other hand, the relation between $W$, $W_1$ and
$W_2$ involves cross terms $C(X,Y)Z$:
\[
W(X,Y)Z = W_1(X_1,Y_1)Z_1 + W_2(X_2,Y_2)Z_2 + C(X,Y)Z.
\]
These can be read off equation (\ref{Intro:curv}). Then $W=0$ implies
$W_1=0$, $W_2=0$ and $C=0$.
We have
\[
C(X_1,Y_1)Z_2 = \frac{1}{2(n+1)}\biggl[ -2\omega(X_1,Y_1)A_2Z_2\biggr]
\]
so $A_2=0$ and interchanging 1 and 2 we see also $A_1=0$. Thus $r_1=0$
and $r_2=0$, and hence $R_1=0$ and $R_2=0$.

\section{Some facts about symmetric symplectic spaces}
Affine symmetric spaces are studied in Loos \cite{Loos}, symplectic
symmetric spaces are studied in Bieliavsky \cite{Biel}.

\begin{definition}\normalfont\rmfamily
A {\bf symmetric symplectic manifold} is a triple $(M,\omega,s)$ where
$M$ is a smooth connected manifold, where $\omega$ is a smooth
symplectic form on $M$ and where $s$ is a smooth map $M\times M \to M$,
$(x,y)\mapsto s_x(y)$, such that:
\begin{itemize}
\item[(i)] for each $x$ in $M$, $s_x$ is an involutive symplectic
diffeomorphism of $(M,\omega)$ (called the symmetry at $x$) and $x$ is
an isolated fixed point of $s_x$,
\item[(ii)] $s_x s_y s_x=s_{s_x(y)}$  for all $x,y$ in $M$.
\end{itemize}
The {\bf transvection group} $G$ of $(M,\omega,s)$ is the group
generated by products of an even number of symmetries.
\end{definition}

We recall below some general facts about symmetric spaces (\cite{Loos},
\cite{Biel}).
 
\Item $(M,\omega,s)$ has a unique connection $\nabla$ such that
$\nabla\omega=0$ and  such that each symmetry $s_x$ is an affine
transformation of $(M,\nabla)$. Observe that ${{s_x}_*}_x=-\Id_{T_xM} $
because $({{s_x}_*}_x)^2=\Id_{T_xM}$ and $x$ is an isolated fixed point
of $s_x$. Since $\omega_x(\nabla_XY,Z) =
\frac{1}{2}(\omega_x(\nabla_XY,Z) +({s_{x}}^* \omega)_x(\nabla_{X} Y,
Z))$, the connection is given by
\begin{equation}\label{sss:connsss}
\omega_x(\nabla_XY,Z) =
 \frac{1}{2}X_x(\omega(Y + s_{x}\cdot Y,Z))
\end{equation}
for $x\in M$, where $X,Y,Z$ are vector fields on $M$ and $(s_{x}\cdot
Y)_y = {s_{x}}_* Y_{s_{x}(y)}$. This connection $\nabla$ has no torsion
and is thus a symplectic connection. The symmetry $s_x$ coincides with
the geodesic symmetry around $x$, since an affinity is determined by its
$1$-jet at a point.

\Item The automorphism group $\Aut = \Aut(M,\omega,s)$ of $(M,\omega,s)$
is the set of symplectic automorphisms $\varphi$ of $(M,\omega)$ such
$\varphi\circ s_x=s_{\varphi(x)}\circ\varphi$, $\forall x\in M$. It is
the intersection of the affine group of $(M,\nabla)$ and the symplectic
diffeomorphism group of $(M,\omega)$. It is thus a Lie group containing
the transvection group so acts transitively on $M$ (since any two points
in $M$ can be joined by a broken geodesic).

Choose a base point $o$ in $M$. Denote by $\widetilde{\sigma}$ the
conjugation by the symmetry $s_o$, it is an involutive automorphism of
$\Aut$.

Let $K'$ denote the stabilizer of $o$ in $\Aut$ and let
$A^{\widetilde{\sigma}}$ (respectively $A_o^{\widetilde{\sigma}}$)
denote the group of fixed points of $\widetilde{\sigma}$ in $\Aut$
(respectively its connected component). Then
$A^{\widetilde{\sigma}}\supseteq K' \supseteq A_o^{\widetilde{\sigma}}$.

Hence, if ${\a}$ (respectively ${\k}'$) is the Lie algebra of $\Aut$
(respectively $K'$) and if $\sigma =\widetilde{\sigma}_{\star Id}$, then
${\k}'$ is the subalgebra of ${\a}$ of fixed points of $\sigma$.

\Item Let  ${\p} = \{X\in{\a}\suchthat \sigma(X)=-X\}$. Then ${\a}=
{\k}' \oplus{\p}$.
 
Denote by $\pi'$ the projection $\Aut\to M$ given by $\pi'(g)=
g\cdot o$. Then ${\pi'_*}_e|_{\p}\colon {\p}\to T_oM$ is a linear
isomorphism which identifies the tangent space $T_oM$ with $\p$.

Denote by $\Exp \colon T_oM\to M$ the exponential map given by the
connection $\nabla$ at the point $o$ and by $\exp$ the exponential map
{}from the Lie algebra $\a$ to the Lie group $\Aut$.

Observe that $s_{\Exp \frac{t}{2}X}s_o$, $X\in T_oM$, is an affinity in
$G$  which realises the parallel transport along $\Exp tX$, since
$s_{\Exp uX*}$ for any $u\in\R$ maps a vector field  which is parallel
along the geodesic $\Exp tX$ to another such parallel vector field.
Hence $s_{\Exp \frac{t}{2}\pi'_{\star e}{X}}s_o=\exp tX$, $\forall X\in
\p$.

It follows that the transvection group $G$, which is stable by
$\widetilde{\sigma}$, is the connected Lie subgroup of $\Aut(M,\omega,s)$
whose Lie algebra is
\begin{equation}
\g=\k\oplus \p  \qquad {\rm{where}}\qquad\k=[\p,\p].
\end{equation}
Indeed, if $G_1$ denotes that subgroup, clearly  by the above
$G_1\subset G$ and the parallel transport along a geodesic $\Exp tX$ is
in $G_1$, but then  any $x\in M$ can be written as $x=g\cdot o$ for
$g\in G_1$ hence $s_xs_o=gs_og^{-1}s_o=g\widetilde{\sigma}(g^{-1})\in
G_1$ and $G\subset G_1$.

Let $K$ denote the stabilizer of $o$ in $G$. Its Lie algebra is $\k$ and
$\k=\{ X\in \g \suchthat \sigma(X)=X\}$. Since the Lie group $G$ acts
effectively on $M$, the representation of $K$ on $T_oM$, $k\mapsto
{k_*}_o$, is faithful so $\k$ acts faithfully on $\p$.

\Item Denote by $\pi$ the projection $\pi \colon G\to M$ where
$\pi(g)=g\cdot o$. Denote by $X^*$ the vector field on $M$ which is the
image under $\pi_*$ of the right invariant vector field on $G$, i.e.
${X^*}_{g\cdot o}=\frac{d}{dt}\exp tX\cdot g\cdot o\vert_{t=0}$.
Observe that $[X^*,Y^*]=-[X,Y]^*$. 
Since $\omega$ is invariant under $G$, formula (\ref{sss:connsss})
yields $\omega_x(\nabla_{Y^*}X^*,Z^*)=\frac{1}{2}
\omega_x([{Y^*},X^*+s_{x}\cdot X^*],Z^*)$
 so $(\nabla_{X^*}Y^*)_{x}=[{X^*},Y^*]+\frac{1}{2}[{Y^*}, 
 X^*+s_{x}\cdot X^*]$.
But $s_{g\cdot o}\cdot X^*=g\cdot s_{o}\cdot g^{-1}\cdot X^* = 
(\Ad g\sigma(\Ad g^{-1}X))^*$
so the connection has the form
\begin{equation}
 (\nabla_{X^*}Y^*)_{g\cdot o}
 =([Y,\Ad g(\Ad g^{-1}X)_\p])^*_{g\cdot o}
\end{equation}
where $Z_\p$ denotes the component in $\p$ of $Z\in \g$ relatively to
the decomposition $\g=\p \oplus \k$ and where $[~,~]$ is the bracket in
${\g}$.
 
Since any $G$-invariant tensor on $M$ is parallel, the curvature tensor
of $(M,\nabla)$ is parallel $(\nabla R=0$) and if $X, Y, Z$ belong to
${\p}$, one has,
\begin{equation}\label{sss:curv}
 R_o(X^*_o,Y^*_o)Z^*_o=-([[X,Y],Z])^*_o.
\end{equation}

\begin{definition} \normalfont\rmfamily
A {\bf symmetric symplectic triple} is a triple $({\g}, \sigma, \Omega)$
where ${\g}$ is a finite dimensional real Lie algebra, $\sigma$ is an
involutive automorphism of ${\g}$ such that if we write ${\g}={\k}\oplus
{\p}$ with $\sigma=\Id_{\k}\oplus -\Id_{\p}$, then
\begin{itemize}
\item $[{\p}, {\p}]= {\k}$;
\item the action of ${\k}$ on ${\p}$ is faithful
\end{itemize}
and where $\Omega$ is a non degenerate skewsymmetric 2-form on $\p$,
invariant by $\k$ under the adjoint action.
\end{definition}

We have seen above that to any connected symmetric symplectic manifold 
$(M,\omega,s)$, when one chooses a base point $o\in M$, 
one associates a symmetric symplectic triple $({\g}, \sigma, \Omega)$
with $\g$  the Lie algebra of its transvection group, with $\sigma$
the differential at the identity of the conjugation by the symmetry $s_o$ 
and with $\Omega=\omega_o$ with the identification between $T_oM$ and $\p$.

Reciprocally, given a symmetric symplectic triple $({\g}, \sigma,
\Omega)$, one builds a simply-connected symmetric symplectic space
$(M,\omega,s)$ with $M=G/K$ where $G$  is the simply-connected Lie group
with Lie algebra $\g$ and $K$ is its connected subgroup with Lie algebra
$\k$, with $\omega$ the $G$-invariant $2$-form on $M$ whose value at
$eK$ is given by $\Omega$ (identifying $T_{eK}M$ and $\p$ via the
differential of the canonical projection $\pi \colon G\to G/K$) and with
symmetries defined by $s_{\pi(g)}\pi(g') = \pi(g\widetilde{\sigma}
(g^{-1}g'))$ where $\widetilde{\sigma}$ is the automorphism of $G$ whose
differential at $e$ is $\sigma$.

\def\temprefb{\ref{Intro:thm2}}
\section{Proof of Theorem \temprefb}

Consider a symmetric symplectic space $(M,\omega,s)$ and assume that its
canonical invariant symplectic connection $\nabla$ has a curvature with
$W=0$.

Since $\nabla R=0$, the Ricci tensor $r$ and  its associated
endomorphism $A$ (where $r(X,Y)=\omega(X,AY)$) are covariantly constant
and hence $A$ commutes with the curvature endomorphisms
\[
AR(X,Y) = R(X,Y)A.
\]
This implies, when we substitute $R$ by its expression in terms of $A$
into equation (\ref{Intro:curv})
\begin{eqnarray*}
-\omega(X,Z)A^2Y + \omega(Y,Z)A^2X = \omega(Y,A^2Z)X - \omega(X, A^2Z)Y.
\end{eqnarray*}
If $Y\ne0$ is arbitrary, $Z=Y$, and we pick $X$ so that $\omega(X,Y)=1$,
then $\omega(Y,A^2Y)=\omega(AY,AY)=0$, so $A^2Y = \lambda_Y Y$ for some 
function $\lambda_Y$. Substituting back into the equation shows that 
$\lambda_Y=\lambda$ is independent of $Y$, and since $A$ is covariant 
constant, $\lambda$ must be constant.

Remark that if $\lambda \ne 0$ then  $r$ is a non-degenerate parallel
symmetric bilinear form so $\nabla$ is its Levi-Civita connection and
$(M,r,s)$ is a pseudo-Riemannian symmetric space.

Let $G$ be the transvection group of our symmetric symplectic space.
Choose a base point $o\in M$ and let $(\g, \sigma, \Omega)$ be the
symmetric triple associated to $(M,\omega, s)$. Let $\g = \k + \p$ be
the decomposition of the Lie algebra of $G$ into the $+1$ and $-1$
eigenspaces of $\sigma$. Then $\Omega(X,Y)=\omega_o(X^*_o,Y^*_o)$ and
with a slight abuse of notations we denote by $R$ the map
$R \colon \p\times\p \to \End(\p)$ so that $(R(X,Y)Z)^*_o =
R_o(X^*_o,Y^*_o)Z^*_o$ and by $A$ the map $A \colon \p \to \p$
so that $(A(X))^*_o=A_o(X^*_o)$. Since $\k$ acts faithfully on $\p$, we
view $\k$ as a subset of $\End(\p)$; by formula (\ref{sss:curv}),
\begin{equation}
\k=\{R(X,Y) \in \End(\p)\suchthat X,Y\in\p\}
\end{equation}
and the brackets on $\g \subset \p\oplus\End(\p)$ are
\begin{equation}\label{thm2:bracket}
[(C,X),(D,Y)]=([C,D]-R(X,Y),CY-DX) 
\end{equation}
where $C,D \in\k\subset\End(\p)$, and $X,Y\in\p$.

Define the $1$-form on $\p$ corresponding to a vector $X \in \p$ by
$\ul{X} = i(X)\Omega$. Formula (\ref{Intro:curv}) giving the curvature
when $W=0$ is equivalent to
\[
R(X,Y) =  k\left(2\Omega(X,Y)A + AY\otimes \ul{X} - AX\otimes
\ul{Y} + X\otimes \ul{AY} - Y\otimes \ul{AX}\right)
\]
where $k = 1/(m+2)$  if $m= \dim M = 2n$. Note that, for a symplectic
symmetric space built from a Lie algebra $\g=\k+\p$ whose bracket of
$\p$ into $\k$ is given by this formula, then the canonical connection
will have curvature given by this formula and hence $W$ will vanish.

Define $B = Y\otimes \ul{X} - X\otimes\ul{Y}$. Clearly $B$ satisfies
$\Omega(U,BV) = \Omega(BU,V)$ and any antisymplectic endomorphism of
$\p$ can be written as a sum of such operators. Then
\[
R(X,Y) =k(\Tr(B)A + AB + BA) 
\]
and, if we put $B'= k(B + \frac12 \Tr(B)I)$, the RHS becomes $C=AB'+B'A$.

\begin{lemma}
For any $\lambda$, 
\[
\k=\{C=AB+BA \suchthat B\in \End(\p) \mathrm{~and~} \Omega(X,BY) =
\Omega(BX,Y)\}.
\]
If $\lambda \ne 0$ then $\k$ is the set of endomorphisms $C\in \End(\p)$
which are infinitesimally symplectic and commute with $A$.
\end{lemma}

\begin{proof}
The first part follows from the considerations above and the fact that
the map $B\mapsto B + \frac12 \Tr(B)$ is a bijection on the space of
antisymplectic endomorphisms of $\p$. $C$ commutes with $A$ since $AC =
\lambda B' + AB'A = CA$. Also $\Omega(X,CY) =  -\Omega(AX,B'Y) +
\Omega(B'X,AY) = -\Omega(B'AX,Y) - \Omega(AB'X,Y) = -\Omega(CX,Y)$.

Conversely, if $\lambda \ne 0$, given $C$ commuting with $A$ and such
that $\Omega(X,CY)= -\Omega(CX,Y)$, let $B=\frac12 \lambda^{-1}AC$; then
\[
BA+AB = \frac12 \lambda^{-1}2\lambda C
= C.
\]
\end{proof}

\subsection{Case $\lambda > 0$}

Write $\lambda = a^2$, $a>0$. Then $\p = V^+ \oplus V^-$ where $V^\pm
=\{ X \in \p \suchthat AX = \pm aX\}$. Let $P^\pm$ be the projection
onto $V^\pm$. Then $A=a(P^+-P^-)$. Clearly
\[\begin{array}{rl}
\omega(V^+,V^+) &= \omega(V^-,V^-) = 0\\
R(V^+,V^+) &= R(V^-,V^-) = 0\\
R(X,Y)&= 2ka(\Omega(X,Y)(P^+ - P^-) - Y\otimes\ul{X} -
X\otimes\ul{Y}),
\end{array}\]
for $X\in V^+$, $Y \in V^-$. It follows that $V^\pm$ are Lagrangian
subspaces of $\p$. Identifying $V^-$ with $(V^+)^*$
via $Y \mapsto \ul{Y}|_{V^+}$ and renaming $V^+$ as $V$, we have
identified $\p$ with $V\oplus V^*$ with its standard symplectic
structure $\Omega(X+\xi, X'+\xi')=-\langle X,\xi'\rangle + \langle
X',\xi \rangle$, and $A$ acts as $+a$ on $V$, $-a$ on $V^*$. With this
notation the curvature has the form
\[
R(X,\xi) = 2ak(-\langle X,\xi \rangle(\Id_{V} - \Id_{V^{*}}) +
\xi\otimes X - X \otimes\xi).
\]

The symplectic centraliser of $A$ can then be identified with
$\End(V)=\gl(V)$, identifying the element in $\End(\p) = \End (V \oplus
V^*)$ given by
\[
\left(\begin{array}{cc}C  &0\\
0 & -{}^tC \end{array}\right)
\]
with the element $C\in\gl(V)$.

So $\k=\gl(V)$ and as a vector space $\g = \gl(V) \oplus V \oplus V^*$
with the brackets 
\begin{eqnarray*}
[(C,X,\xi),(C',X',\xi')] &=& ([C,C']+2ka(\langle X,\xi'\rangle - \langle
X',\xi \rangle)I \\&&\qquad +2ka X\otimes\xi' - 2ka X'\otimes\xi,\\
&&\qquad\qquad CX' - C'X, -{}^tC\xi' + {}^tC'\xi).
\end{eqnarray*}
The map $j \colon \g \to \esl(V\oplus\R)$ given by
\[
j(C,X,\xi) = \left(\begin{array}{cc}
C - 2k \Tr(C) I & sX\\
s\, {}^t \xi & - 2k\Tr(C)
\end{array}\right)
\]
has the brackets above provided $s^2=2ka$.

Thus when $\lambda > 0$, $M  = G/K$ where $G = SL(n+1,\R)$, $K =
GL(n,\R)$. The involution $\sigma$ is given by
\[
\sigma\left(\begin{array}{cc}C&v\\ \xi&-\Tr(C)\end{array}\right) = 
\left(\begin{array}{cc}C&-v\\-\xi&-\Tr(C)\end{array}\right)
\]
and, writing $(X,\xi)$ for $\left(\begin{array}{cc}0&X\\
\xi&0\end{array}\right)$, the symplectic form  is given by
\[
\Omega\left((X,\xi),(X',\xi')\right)= 
-\langle X,\xi'\rangle + \langle X',\xi\rangle.
\]
The curvature of the canonical connection on 
this symplectic symmetric space at the base point $eK$ is
\begin{eqnarray*}
&~&R((X,\xi),(X',\xi'))(X'',\xi'')
= (X''(\langle X',\xi\rangle
-\langle X,\xi'\rangle)-X\langle X'',\xi'\rangle\\ 
&&\quad+X'\langle X'',\xi\rangle\, ,\,\xi'\langle
X,\xi''\rangle-\xi\langle X',\xi''\rangle -\xi''(\langle X',\xi\rangle
-\langle X,\xi'\rangle))
\end{eqnarray*}
\[
r((X,\xi),(X',\xi'))=(n+1)(\langle X,\xi'\rangle +\langle X',\xi\rangle)
\]
\[
A(x,\xi)=(n+1)(x,-\xi)
\]
and formula (\ref{Intro:curv}) holds so $R$ is of Ricci-type.

\subsection{Case $\lambda < 0$}

We write $\lambda = -b^2$ where $b<0$. If we put $J= b^{-1} A$ then $J$
defines a complex structure on the vector space $\p$. We write $V$ for
$\p$ viewed as an $n$-dimensional complex vector space. $V$ has a
(pseudo-)Hermitean structure given by
\[
\langle X,Y \rangle = \Omega(X,JY) + i \Omega(X,Y)
\]
which is $\C$-linear in the second variable. The infinitesimally
symplectic transformations which commute with $A$, or equivalently $J$,
are the complex linear transformations of $V$ which are skew-Hermitean
with respect to this Hermitean structure. Thus $\k$ is the (pseudo-)
unitary Lie algebra $\u(V,\langle\,,\,\rangle)$.

The curvature has the form
\[
R(X,Y) = kb \left( 2\Omega(X,Y)J + Y \otimes\langle X, \,.\, \rangle 
- X \otimes\langle Y, \,.\, \rangle \right).
\]

Then $\g = \u(V,\langle\,,\,\rangle) \oplus V$ with bracket
\begin{eqnarray*}
[(C,X),(C',X')] &=& ([C,C']+kb(X\otimes\langle X', \,.\,\rangle - 
X'\otimes\langle X, \,.\,\rangle\\
&&\qquad -2\Omega(X,X')J),CX'-C'X).
\end{eqnarray*}
and $\g$ can be identified with $\su(V\oplus\C, \llangle\,,\,\rrangle)$
via
\[
j(C,X)=\left( \begin{array}{cc}
C - 2k\Tr(C) I & sX\\
-\ol{s} \langle X, \,.\, \rangle & -2k \Tr(C)
\end{array} \right)
\]
with
\[
\llangle (v,r), (w,t) \rrangle = \langle v,w \rangle + \ol{r}t
\]
provided
\[
s\ol{s} = -kb.
\]
Hence when $\lambda<0$ then $M=G/K$ with $\g=\su(p+1,q)$, $p+q=n$,
$\k= \u(p,q)$,
\[
\sigma\left( 
\begin{array}{cc}C & v\\
-\langle v, \,.\, \rangle & -\Tr(C)\end{array} 
\right) = \left( 
\begin{array}{cc}C & -v\\\langle v, \,.\, \rangle & -\Tr(C)\end{array} 
\right)
\]
where $\langle v,w \rangle=\sum_{i=1}^p\overline{v}^i w^i-\sum_{j=p+1}^n\overline{v}^j w^j$ and
\[
\Omega(v,w) = \Im \langle v,w \rangle.
\]
The curvature of the canonical connection on this symmetric symplectic
space at $eK$ is
\[
R(v,w)z
= v\langle w,z\rangle - w\langle v,z\rangle + z(-\langle v,w\rangle + 
\langle w,v\rangle)
\]
\[
r(v,z)=-2(n+1)\Re \langle v,z\rangle 
\]
\[
A(v)=-2(n+1)iv
\]
and formula (\ref{Intro:curv}) holds so $R$ is of Ricci-type.

\subsection{Case $\lambda=0$}

In this case $A$ is nilpotent since $A^2=0$. Let $Z = \Image A$ and
$\widetilde{Z} = \Ker A$. Then $Z \subset \widetilde{Z}$, and $Z$ and
$\widetilde{Z}$ are symplectic orthogonals of each other. If $V$ denotes
a complement for $Z$ in $\widetilde{Z}$, then the restriction of
$\Omega$ to $V$ is non-degenerate. $Z$ is contained in the
$\Omega$-orthogonal of $V$; let $Z'$ be a complement so that $V^\perp =
Z \oplus Z'$. $V^\perp$ is a symplectic subspace and $Z$ is maximal
isotropic so we can also suppose that $Z'$ is maximal isotropic.
$\Omega$ gives a duality of $Z$ with $Z'$.

In other words, we have written $\p$ as $Z \oplus Z^* \oplus V$ where $Z
\oplus Z^*$ has its standard symplectic structure and $V$ is a
symplectic vector space. $A$ is non-zero only on $Z^*$ and maps it
isomorphically onto $Z$, and as such it is symmetric. In block form, the
symplectic structure $\Omega$ is given by
\[
\left( 
\begin{array}{ccc}
0 & -I & 0\\
I & 0 & 0\\
0 & 0 & J'\end{array} 
\right)
\]
and $A$ by
\[
\left( 
\begin{array}{ccc}
0 & A' & 0\\
0 & 0 & 0\\
0 & 0 & 0\end{array} 
\right)
\]
where $A'$, by a suitable choice of basis is diagonal with $\pm1$ on the
diagonal. An easy calculation shows that matrices of the form $AB+BA$ 
with $\Omega(X,BY)=\Omega(BX,Y)$  have the form
\[
\left( 
\begin{array}{ccc}
K & L & -{}^tMJ'\\
0 & -{}^tK & 0\\
0 & M & 0\end{array} 
\right)
\]
where ${}^tKA' + A' K=0$, ${}^tL=L$. The matrices with $K=0$ form an
ideal which is 2-step nilpotent (abelian when $\rank A=n=\frac{1}{2}
\dim M$) and the matrices with $L=M=0$ a subalgebra isomorphic to
$\so(p,q)$, where $p+q=r=\rank A$, $p$ the number of $+$'s and $q$ the
number of $-$'s in $A'$ (hence $(p,q)$ is the signature of the non
degenerate symmetric bilinear form naturally induced on $\p/\Ker A$ by
the Ricci tensor $\Omega(X,AY)$).

The bracket of $\p$ into $\k$ is given, using formulas 
(\ref{thm2:bracket}) and (\ref{Intro:curv})  by
\[
\left[ \left(\begin{array}{c}u\\v\\w\end{array}\right), 
\left(\begin{array}{c}u'\\v'\\w'\end{array}\right)\right]
=-k\left( 
\begin{array}{ccc}
{\tilde{K}}=A'(v'{}^tv - v{}^tv') &{\tilde{L}}& -{}t{\tilde{M}}J')\\
0 & -{}^t{\tilde{K}} & 0\\
0 & {\tilde{M}} & 0\end{array} 
\right)
\]
where ${\tilde{L}}=A'B+{}^tBA'+2(\Tr B+{}^twJ'w')A'$ with
$B=v{}^tu'-v'{}^tu$ and ${\tilde{M}}=-{}^t(A'(v'{}^tw - v{}^tw'))$.

Then $\g=\k+\p=\{(K,L,M,u,v,w) \suchthat K\in \so(p,q), L\in 
\Mat(r\times r,\R), {}^tL=L, M\in\Mat(2n-2r\times r,\R), u \in Z={\R}^r,
v\in Z^*, w\in W={\R}^{2n-2r} \}$.
The brackets are given, with obvious notations, by
\[
[(K,L,M),(K',L',M')]=([K,K'],L'',-M{}^tK'+M'{}^tK)
\]
where $L''=KL'-L{}^tK' -K'L+L'{}^tK-{}^tMJ'M'+{}^tM'J'M$,
\[
[(K,L,M),(u,v,w)]=(Ku+Lv-{}^tMJ'w,-{}^tKv,Mv),
\]
\[
[(u,v,w),(u',v',w')] = (-kA'(v'{}^tv - v{}^tv'),-k{\tilde{L}},
k\,{}^t(A'(v'{}^tw - v{}^tw')))
\]
where ${\tilde{L}}$ is defined as above.

We can combine $\so(p,q)$ with $Z^*$ to give $\so(p,q+1)$ via
\[
(K,v) \mapsto \left(\begin{array}{cc} K &-k^{1/2} A'v\\
-k^{1/2}\, {}^tv &0\end{array}\right).
\]
The subset $\r=\{ (0,L,M,u,0,w) \in \g\}$ is a 2-step nilpotent ideal of
$\g$ (abelian when $r=n$ i.e. when the rank of the Ricci tensor is half
the dimension of the manifold). Hence, when $p+q=r>1$, $\r$ is the
radical of $\g$ and the semisimple Levi factor of $\g$ is isomorphic to
$\so(p,q+1)$.

\section{Some corollaries}

\begin{cor}
Let $(M_i,\omega_i,s_i)$, $i=1,2$ be symmetric symplectic spaces of the
same dimension $2n$ with $W_i=0$ with semisimple transvection groups $G_i$. 
Then $G_1^{\C}=G_2^{\C}$.
\end{cor}

\begin{proof}
$SL(n+1,\R)$ and $SU(p+1,q)$ both have $SL(n+1,\C)$ as complexification.
\end{proof}

\begin{cor}
Let $(M,\omega,s)$ be a compact, simply-connected symmetric symplectic space
of dimension $2n$ such that $W=0$ then $(M,\omega,s)$ is $\P_n(\C)$.
\end{cor}

\begin{proof}
This follows immediately from the list in Theorem \ref{Intro:thm2}. The only
case where $G/K$ is compact is when $G=SU(n+1)$ and $K=U(n)$.
\end{proof}

In dimension 4 we have the following list of possibilities (up to coverings)
for $M$:

\begin{itemize}
\item $SL(3,\R)/GL(2,\R)$;
\item $SU(1,2)/U(2)$;
\item $SU(2,1)/U(1,1)$;
\item $SU(3)/U(2)$;
\item $\lambda=0$ cases corresponding to:
\begin{itemize}
\item[{$\circ$}] $\rank A=1$, $p=0$ or $p=1$;
\item[{$\circ$}] $\rank A=2$, $p=0$, $p=1$ or $p=2$.
\end{itemize}
\end{itemize}

\end{document}